\documentclass[12pt]{article}

\usepackage{a4wide}

\usepackage[intlimits]{amsmath}
\usepackage{amssymb,amsthm}

\usepackage[cp1251]{inputenc}
\usepackage[english]{babel}

\newtheorem{thm}{Theorem}[section]

\newtheorem{lem}[thm]{Lemma}

\theoremstyle{definition}

\theoremstyle{remark}
\newtheorem{rem}[thm]{Remark}
\numberwithin{equation}{section}

\begin{document}\Large

\title{Asymptotic expansion of Markov random evolution}

\author{I.V.Samoilenko\\
Institute of Mathematics,\\ Ukrainian National Academy of
Sciences,\\ 3 Tereshchenkivs'ka, Kyiv, 01601, Ukraine\\
isamoil@imath.kiev.ua}

\maketitle

\abstract{Is studied asymptotic expansion for solution of
singularly perturbed equation for Markov random evolution in
$\mathbb{R}^d$.
The views of regular and singular parts of solution are found. }

\section{\bf Introduction}

A Markov random evolution (MRE) is created by a solution of the
evolutionary equation in Euclidean space $\mathbb{R}^d, d\geq1$
$$du^{\varepsilon}(t)/dt=v(u^{\varepsilon}(t);{\ae}(t/\varepsilon))$$
with the ergodic Markov switching process ${\ae}(t), t\geq0$ on the
standard (Polish) phase-space $(E,\mathcal{E})$ by the operator
$Q(x,B), x\in E, B\in \mathcal{E}$ that defines transition
probabilities of a Markov chain $\ae_n, n\geq 0$
$$Q(x,B)=P\{\ae_{n+1}\in B|\ae_n=x\}.$$

The operator of transition probabilities $Q$ is defined by
$$Qf(x)=\int_EQ(x,dy)f(y), x\in E, \eqno(1)$$ for any bounded measurable real
valued $f$ defined on $E$.

We will see later that the equation for the regular and the
singular parts of a random evolution are defined by the generator
(1) of a uniformly argodic Markov switching process. The Banach
space $\mathcal{B}(E)$ is splitted onto the two subspaces {\rm
\cite{KorTur}}:
$$\mathcal{B}(E)=N_Q\bigoplus R_Q,$$ where
$N_Q:=\{\varphi:Q\varphi=0\}$ is the null-space of $Q$, and
$R_Q:=\{\psi:Q\varphi=\psi\}$ is the range of $Q$.

We define the projector $\Pi:$ $N_Q:=\Pi \mathcal{B}(E),
R_Q:=(I-\Pi)\mathcal{B}(E);
\Pi\varphi(x):=\widehat{\varphi}\mathbf{1},
\widehat{\varphi}:=\int_E \varphi(x)\pi(dx),$ where the stationary
distribution $\pi(B), B\in \mathcal{E}$ of the Markov process
${\ae}(t),t\geq0$ satisfies the relations {\rm \cite{KoKo}}
$$\pi(dx)=\rho(dx)m_1(x)/\widehat{m},$$
$$\widehat{m}=\int_E m_1(x)\rho(dx).$$ $\rho(B), B\in \mathcal{E}$ is
the stationary distribution of the Markov chain  ${\ae}_n, n\geq
 0,$ given by the equation $$\rho(B)=\int_E Q(x,B)\rho(dx), \rho(E)=1.$$

Let us consider the Banach space $\mathcal{B}(\mathbb{R}^d)$ of
real-valued test-functions $\varphi(u), u\in \mathbb{R}^d$ which are
bounded with all their derivatives equipped with $\sup$-norm
$$||\varphi||:=\sup\limits_{u\in
\mathbb{R}^d}|\varphi(u)|<C_{\varphi}.$$

The random evolution in $\mathcal{B}(\mathbb{R}^d)$ is given by the
relation
$$\Phi_t^{\varepsilon}(u,x):=E[\varphi(u^{\varepsilon}(t))|u^{\varepsilon}(0)=u, {\ae}^{\varepsilon}(0)=x]. \eqno(2)$$

The asymptotic behavior of MRE (2) as $\varepsilon\to 0$ is
investigated under the assumption of uniformly ergodicity of the
Markov switching process ${\ae}(t)$ described above and under the
assumption of the existence of a global solution of the
deterministic equations
$$du_x(t)/dt=v(u_x(t);x), x\in E.$$

Let us consider the deterministic evolution $$\Phi_x(t,u)=\varphi
(u_x(t)), u_x(0)=u.$$ It generates a corresponding semigroup
$$\mathbb{V}_t(x)\varphi(u):=\varphi (u_x(t)), u_x(0)=u,$$ and its
generator has the form:
$$\mathbb{V}(x)\varphi(u)=v(u;x)\varphi '(u):=\sum_{k=1}^dv_k(u;x)\varphi_k '(u),$$
$$\varphi_k '(u):=\partial\varphi(u)/\partial u_k, \varphi(u)\in C^{\infty}(\mathbb{R}^d).$$

By the average principle {\rm \cite{Skor}} the weak convergence
$$u^{\varepsilon}(t)\Rightarrow \widehat{u} (t), \varepsilon\to 0 \eqno(3)$$
takes place. The average limit evolution $\widehat{u} (t), t\geq 0$
is defined by a solution of the average equation $$d\widehat{u}
(t)/dt=\widehat{v}(\widehat{u} (t)).$$

The average velocity $\widehat{v}(u), u\in \mathbb{R}^d$ is defined
by $$\widehat{v}(u)=\int_Ev(u;x)\pi(dx)$$ (i.e. by the average of
the initial velocity $v(u;x)$ over the stationary distribution
$\pi(B), B\in \mathcal{E}$).

The rate of convergence in (3) can be investigated in two
directions: \\ i) asymptotic analysis of the fluctuations
$$\zeta^{\varepsilon}(t)=u^{\varepsilon}(t)-\widehat{u}(t); \eqno(4)$$
\\ ii) asymptotic analysis of the average deterministic evolution
(2).

The asymptotic analysis of fluctuations (4) leads to the diffusion
approximation of the random evolution {\rm \cite{KoLim, Skor}}.

The asymptotic analysis of evolution (2) is realized in what
follows by constructing the asymptotic expansion in power of the
small parameter series $\varepsilon\to 0 (\varepsilon>0)$ in the
following form ($\tau=t/\varepsilon$):
$$\Phi_t^{\varepsilon}(u,x)=u^{(0)}(t)+\sum_{k=1}^{\infty}\varepsilon^k[u^{(k)}(t)+w^{(k)}(\tau)]. \eqno(5)$$

The asymptotic expansion (5) contains two parts: \\ i) the regular
term $u^{\varepsilon}(t):=u^{(0)}(t)+\sum_{k\geq
1}^{\infty}\varepsilon^k u^{(k)}(t),$ \\ ii)the singular term
(boundary layer) $w^{\varepsilon}(\tau):=\sum_{k\geq
1}^{\infty}\varepsilon^k w^{(k)}(\tau), \tau=t/\varepsilon.$

In addition the initial condition:
$$u^{(0)}(0)=\varphi(u)\mathbf{1}$$
has to be valid for any $x\in E, u\in \mathbb{R}^d$.


It's well-known (see, e.g. {\rm \cite{Pin}}), that the evolution,
determined
 by a test-function $\varphi(u)\in C^{\infty}(R^d)$ (here $\varphi(u)$
is integrable on $\mathbb{R}^d$
): satisfy the system
of Kolmogorov backward differential equations:
$$\begin{array}{c}
    \frac{\partial}{\partial
t}\Phi_t^{\varepsilon}(u,x)=[\varepsilon^{-1}Q+\mathbb{V}]\Phi_t^{\varepsilon}(u,x), \\
   \Phi_0^{\varepsilon}(u,x)=\varphi(u).
  \end{array}
 \eqno(6)$$

Asymptotic expansions with 'boundary layers' were studied by many
authors (see {\rm \cite{Kor, KorBor, VasBut}}). In particular,
functionals of Markov and semi-Markov processes are investigated
from this point of view in {\rm \cite{KorPenTur, Sam, Tad}}.

In this work we study system (6) with the first order singularity.
To find asymptotic expansion of the solution of (2) we use the
method proposed in {\rm \cite{KorBor, VasBut}}. The solution
consists of two parts - regular terms and singular terms - which
are determined by different equations. Asymptotic expansion lets
not only determine the terms of asymptotic, but to see the
velocity of convergence in hydrodynamic limit.

Besides, when studying this problem, we improved the algorithm of
asymptotic expansion. Partially, the initial conditions for the
regular terms of asymptotic are determined without the use of
singular terms, i.e. the regular part of the solution may be found
by a separate recursive algorithm; scalar part of the regular term
is found and without the use of singular terms. These and other
improves of the algorithm are pointed later.

\section{\bf Asymptotic expansion of the solution}

Let $P(t)=e^{Qt}=\{p_{ij}(t); i,j\in E\}.$ Put
$\pi_j=\lim\limits_{t\to\infty}p_{ij}(t)$ and
$-R_0=\{\int_0^{\infty}(p_{ij}(t)-\pi_j)dt; i,j\in E\}=\{r_{ij};
i,j\in E\}.$

Let $\Pi$ be a projecting operator on the null-space $N_Q$ of the
operator $Q$. For any vector $g$ we have $\Pi
g=\widehat{g}\mathbf{1},$ where $\widehat{g}=(g,\pi),
\mathbf{1}=(1,\ldots,1).$ Then for the operator $Q$ the following
correlations are true (see {\rm \cite{KorTur}}, chapter 3)
$$\Pi
Q\Pi=0,$$
$$QR_0=R_0Q=\Pi-I.$$

We put:
$$exp_0(Qt):=e^{Qt}-\Pi,$$

\begin{thm} The solution of equation (6) with initial condition
$\Phi^{\varepsilon}_0(u,x)=\varphi(u)$, where $\varphi(u)\in
C^{\infty}(R^d)$ and integrable on $\mathbb{R}^d$ has asymptotic
expansion
$$\Phi^{\varepsilon}_t(u,x)={u}^{(0)}(t)+\sum_{n=1}^{\infty}\varepsilon^n\left({u}^{(n)}(t)+
{w}^{(n)}\left(t/\varepsilon\right)\right).\eqno(7)$$

Regular terms of the expansion are: ${u}^{(0)}(x,t)$ - the
solution of equation
$$\frac{\partial}{\partial t}{u}^{(0)} (t)-\Pi
\mathbb{V}\Pi {u}^{(0)}(t)=0 \eqno(8)$$ with initial condition
${u}^{(0)}(0)=\varphi(u),$
$${u}^{(1)}(t)=R_0\left[\frac{d}{dt}{u}^{(0)}(t)-\mathbb{V}{u}^{(0)}(t)\right]+c^{(1)}(t):=R_0\mathbb{L}u^{(0)}(t)+c^{(1)}(t),$$ for $k\geq 2:$
$${u}^{(k)}(t)=R_0\mathbb{L}{u}^{(k-1)}(t)+c^{(k)}(t)$$ where
$c^{(k)}(t)\in N_{Q},$
$$c^{(k)}(u,t)=c^{(k)}(V^{-1}(t+V(u)),0)+\int\frac{\mathcal{L}_k(V^{-1}(t+V(u)),0)}{V^{-1}(t+V(u))}du-\int\frac{\mathcal{L}_k(u,t)}{v(u)}du,
k>0,$$ here $V(u)=\int\frac{du}{v(u)},$ $V^{-1}(w)$ is the
backward function for $V(u)$,
$$\mathcal{L}_k(u,t)=\sum_{i=0}^{k-1}\sum_{n=1}^{k-i}(-1)^k(k-i-n+1)\Pi\mathbb{V}\mathbb{R}_0\mathbb{V}^n\Pi\frac{d^{k-i-n}}{dt^{k-i-n}}c^{(i)}(t).$$


The singular terms of the expansion have the
view:$${w}^{(1)}(\tau)=exp_0({Q\tau})\mathbb{V} \varphi(u),$$ for
$k>1:$
$${w}^{(k)}(\tau)=exp_0(Q\tau){w}^{(k)}(0)+\int_0^{\tau}exp_0(Q(\tau-s))\mathbb{V}{w}^{(k-1)}(s)ds-$$ $$\Pi
\int_{\tau}^{\infty}\mathbb{V}{w}^{(k-1)}(s)ds.$$ Initial
conditions:
$$c^{(0)}(0)=\varphi(u),$$
$${w}^{(1)}(0)=-R_0\mathbb{L}\varphi(u),$$
$$c^{(1)}(0)=0,$$
for $k> 1:$
$${w}^{(k)}(0)=-R_0\mathbb{L}{u}^{(k-1)}(0),$$
$$c^{(k)}(0)=\mathbb{V}\widetilde{w}^{(k-1)}(0),$$
where $ \widetilde{w}^{(1)}(0)=-R_0\mathbb{V} \varphi(u),$
$$\widetilde{w}^{(k)}(0)=R_0\mathbb{L}{u}^{(k-1)}(0)+R_0\mathbb{V}\widetilde{w}^{(k-1)}(0)+$$
$$\Pi
\mathbb{V}(\widetilde{w}^{(k-1)}(\lambda))'_{\lambda}|_{\lambda=0},$$
$$(\widetilde{w}^{(k)}(\lambda))'_{\lambda}|_{\lambda=0}=R^2_0\mathbb{L}{u}^{(k-1)}(0)+R_0^2Q_{1}\widetilde{w}^{(k-1)}(0)+$$
$$
R_0\mathbb{V}(\widetilde{w}^{(k-1)}(\lambda))'_{\lambda}|_{\lambda=0}.$$

\end{thm}

\begin{rem} The initial conditions for the regular terms of
asymptotic are determined without the use of singular terms, i.e.
the regular part of the solution may be found by a separate
recursive algorithm (comp. with {\rm \cite{KorPenTur}}).
\end{rem}

{\it Proof of Theorem 2.1:} Let us substitute the solution
$\Phi_t^{\varepsilon}(u,x)$ in the form (7) to the equation (6)
and equal the terms at $\varepsilon$ degrees. We'll have the
system for the regular terms of asymptotic:$$
  \begin{cases}
    Q{u}^{(0)}=0 \\
    Q{u}^{(k)}=\frac{d}{dt}{u}^{(k-1)}-\mathbb{V}{u}^{(k-1)}:=\mathbb{L}{u}^{(k-1)}, k\geq 1
  \end{cases}\eqno(9)
$$   and for the singular terms $$\frac{dw^{\varepsilon}}{dt}=\frac{dw^{\varepsilon}}{d\tau}
\frac{d\tau}{dt}=\varepsilon^{-1}\frac{dw^{\varepsilon}}{d\tau}=(\varepsilon^{-1}Q+\mathbb{V})w^{\varepsilon}.$$

Thus, from $\frac{dw^{\varepsilon}}{d\tau}=(Q+\varepsilon
\mathbb{V})w^{\varepsilon}$ we obtain:
$$
  \begin{cases}
    \frac{d}{d
    \tau}{w}^{(1)}=Q{w}^{(1)} \\
    \frac{d}{d
    \tau}{w}^{(k)}-Q{w}^{(k)}=\mathbb{V}{w}^{(k-1)},
    k>1.
  \end{cases} \eqno(10)
$$

From (9) we have: ${u}^{(0)}(t)\in N_{Q}.$

The solvability condition for ${u}^{(1)}(t)$ has the view: $$\Pi
Q\Pi {u}^{(1)}(t)=0= \frac{\partial}{\partial t}{u}^{(0)}(t)-\Pi
\mathbb{V}\Pi {u}^{(0)}(t).$$

So, we have equation (8) for ${u}^{(0)}(t)$.


For ${u}^{(1)}(t)$ we have:
$${u}^{(1)}(t)=R_0\mathbb{L}{u}^{(0)}(t)+c^{(1)}(t).$$

Using the second equation from (9) we obtain:
$${u}^{(k)}(t)=R_0\mathbb{L}{u}^{(k-1)}(t)+c^{(k)}(t),$$ where $c^{(k)}(t)\in
N_{Q}.$

To find $c^{(k)}(t)$ we'll use the fact that ${u}^{(0)}(t)\in
N_{Q}.$ Let us put $c^{(0)}(t)={u}^{(0)}(t)$.

For the equation $$Q{u}^{(2)}(t)=\frac{\partial}{\partial
t}u^{(1)}(t)-\mathbb{V}u^{(1)}(t)=\frac{d}{dt}{R}_0\left[\frac{d}{dt}c^{(0)}(t)-\mathbb{V}c^{(0)}(t)\right]
+\frac{d}{dt}c^{(1)}(t)-$$ $$
\mathbb{V}{R}_0\left[\frac{d}{dt}c^{(0)}(t)-\mathbb{V}c^{(0)}(t)\right]-\mathbb{V}c^{(1)}(t)$$
we use the solvability condition $$\Pi Q\Pi
{u}^{(2)}(t)=0=\frac{d}{dt}c^{(1)}(t)-\mathbb{V}c^{(1)}(t)+\Pi{R}_0\Pi\frac{d^2}{dt^2}c^{(0)}(t)-
\Pi{R}_0\mathbb{V}\Pi\frac{d}{dt}c^{(0)}(t)-$$ $$
\Pi\mathbb{V}{R}_0\Pi\frac{d}{dt}c^{(0)}(t)+
\Pi\mathbb{V}{R}_0\mathbb{V}\Pi c^{(0)}(t).$$

We find:
$$\frac{d}{dt}c^{(1)}(t)-\mathbb{V}c^{(1)}(t)=-\Pi\mathbb{V}{R}_0\mathbb{V}\Pi c^{(0)}(t).$$

By induction:
$$\frac{d}{dt}c^{(k)}(t)-\mathbb{V}c^{(k)}(t)=\sum_{i=0}^{k-1}\sum_{n=1}^{k-i}(-1)^k(k-i-n+1)\Pi\mathbb{V}
{R}_0\mathbb{V}^n\Pi\frac{d^{k-i-n}}{dt^{k-i-n}}c^{(i)}(t), k>0.$$

So, we have the following equation for $c^{(k)}(u,t)$:
$$\frac{d}{dt}c^{(k)}(u,t)-v(u)\frac{d}{du}c^{(k)}(t)=\mathcal{L}_k(u,t),$$
here
$$\mathcal{L}_k(u,t)=\sum_{i=0}^{k-1}\sum_{n=1}^{k-i}(-1)^k(k-i-n+1)\Pi\mathbb{V}\mathbb{R}_0\mathbb{V}^n\Pi\frac{d^{k-i-n}}{dt^{k-i-n}}c^{(i)}(t).$$

To find a solution we should write down a system
$$\frac{dt}{1}=-\frac{du}{v(u)}=\frac{dc^{(k)}}{\mathcal{L}_k(u,t)}.$$

The independent integrals of this system are:
$$t+\int\frac{du}{v(u)}=C_1,$$
$$c^{(k)}(u,t)+\int\frac{\mathcal{L}_k(u,t)}{v(u)}du=C_2.$$

As soon as $c^{(k)}(u,t)$ is only in one of the first integrals,
we may present the solution in the form:
$$c^{(k)}(u,t)=f_k\left(t+\int\frac{du}{v(u)}\right)-\int\frac{\mathcal{L}_k(u,t)}{v(u)}du,
k>0,$$ where $f_k$ is any differentiable function. Using initial
condition for $c^{(k)}(u,t)$ we find a condition for $f_k$:
$$f_k\left(\int\frac{du}{v(u)}\right)=c^{(k)}(u,0)+\int\frac{\mathcal{L}_k(u,0)}{v(u)}du.$$

We may put now $V(u)=\int\frac{du}{v(u)}$ and make a change of
variables $w=V(u)$. So, $u=V^{-1}(w)$ and we have:
$$f_k(w)=c^{(k)}(V^{-1}(w),0)+\int\frac{\mathcal{L}_k(V^{-1}(w),0)}{v(V^{-1}(w))}\frac{dw}{w}.$$

Thus, we obtain
$$c^{(k)}(u,t)=c^{(k)}(V^{-1}(t+V(u)),0)+\int\frac{\mathcal{L}_k(V^{-1}(t+V(u)),0)}{V^{-1}(t+V(u))}du-\int\frac{\mathcal{L}_k(u,t)}{v(u)}du,
k>0.$$

Initial conditions for $c^{(k)}(u,0)$ are found later through
Laplace transform for the singular terms of asymptotic.


For the singular terms we have from (10):
$${w}^{(1)}(\tau)=exp_0(Q\tau){w}^{(1)}(0).$$

Here we should note that the ordinary solution
${w}^{(1)}(\tau)=exp(Q\tau){w}^{(1)}(0)$ is corrected by the term
$-\Pi{w}^{(1)}(0)$ in order
 to receive the following $\lim\limits_{\tau\to
 \infty}{w}^{(1)}(\tau)=0$.
We choose this limit to be equal 0 for all singular terms, that
may done due to uniform ergodicity of switching Markovian process.

The following statements are made using a method proposed in {\rm
\cite{Kor}}. For the second equation of the system the
corresponding solution should be
$${w}^{(k)}(\tau)=exp_0(Q\tau){w}^{(k)}(0)+
\int_0^{\tau}exp_0(Q(\tau-s))\mathbb{V}{w}^{(k-1)}(s)ds,$$ where the
homogenous part has the following solution
$${w}^{(k)}(\tau)=exp_0(Q\tau){w}^{(k)}(0).$$

But here we should again correct the solution, in order
 to receive the limit
 $\lim\limits_{\tau\to
 \infty}{w}^{(k)}(\tau)=0$, by the term $-
\Pi\int_{\tau}^{\infty}\mathbb{V} {w}^{(k-1)}(s)ds.$

And so the solution is:
$${w}^{(k)}(\tau)=exp_0(Q\tau){w}^{(k)}(0)+
\int_0^{\tau}exp_0(Q(\tau-s))\mathbb{V}{w}^{(k-1)}(s)ds-$$
$$\Pi\int_{\tau}^{\infty}\mathbb{V} {w}^{(k-1)}(s)ds.$$

We should finally find the initial conditions for the regular and
singular terms.

We put $c^{(0)}(t)={u}^{(0)}(t)$, so
$c^{(0)}(0)={u}^{(0)}(0)=\varphi (u).$

From the initial condition for the solution
${u}^{\varepsilon}(0)={u}^{(0)}(0)=\varphi (u)$, we have to
determine ${u}^{(k)}(0)+{w}^{(k)}(0)=0, k\geq 1.$ Let us rewrite
this equation for the null-space $N_Q$ of matrix
$Q$:$$\Pi{u}^{(k)}(0)+\Pi{w}^{(k)}(0)=0, k\geq 1, \eqno(11)$$ and
the space of values $R_Q$:
$$(I-\Pi){u}^{(k)}(0)+(I-\Pi){w}^{(k)}(0)=0, k\geq 1. \eqno(12)$$

So, for $k=1$ we obtain:
$${u}^{(1)}(0)=R_0\mathbb{L}{u}^{(0)}(0)+c^{(1)}(0)=$$
$$=(I-\Pi)R_0\mathbb{L}\varphi(u)+\Pi c^{(1)}(0),$$
$${w}^{(1)}(0)=(I-\Pi){w}^{(1)}(0).$$

Thus, $c^{(1)}(0)=0, {w}^{(1)}(0)=-R_0\mathbb{L}\varphi(u).$

By analogy, for $k>1$:
$${u}^{(k)}(0)=R_0\mathbb{L}{u}^{(k-1)}(0)+c^{(k)}(0)=$$ $$=(I-\Pi)R_0\mathbb{L}{u}^{(k-1)}(0)+\Pi c^{(k)}(0),$$
$${w}^{(k)}(0)=(I-\Pi){w}^{(k)}(0)-\Pi\int_0^{\infty}\mathbb{V}{w}^{(k-1)}(s)ds.$$

Functions ${w}^{(k-1)}(s),{u}^{(k-1)}(0)$ are known from the
previous steps of induction. So, we've found $\Pi{w}^{(k)}(0)$ in
(11) and $(I-\Pi){u}^{(k)}(0)$ in (12).

Now we may use the correlations (11), (12) to find the unknown
initial conditions:
$$c^{(k)}(0)=\int_0^{\infty}\mathbb{V}{w}^{(k-1)}(s)ds,$$
$${w}^{(k)}(0)=-R_0\mathbb{L}{u}^{(k-1)}(0).$$

In {\rm \cite{KorPenTur}} an analogical correlation was found for
$c^{(k)}(0)$. To find $c^{(k)}(0)$ explicitly and without the use
of singular terms we'll find Laplace transform for the singular
term. The following lemma is true.

\begin{lem} Laplace transform for the singular term of asymptotic
expansion
$$\widetilde{w}^{(k)}(\lambda)=\int_0^{\infty}e^{-\lambda s}
{w}^{(k)}(s)ds$$ has the view:

$$\widetilde{w}^{(1)}(\lambda)=(\lambda-\Pi+(R_0+\Pi)^{-1})^{-1}[-R_0\mathbb{V}
\varphi(u)],$$

$$\widetilde{w}^{(k)}(\lambda)=(\lambda-\Pi+(R_0+\Pi)^{-1})^{-1}\mathbb{L}{u}^{(k-1)}(0)+$$
$$(\lambda-\Pi+(R_0+\Pi)^{-1})^{-1}\mathbb{V}
\widetilde{w}^{(k-1)}(\lambda)+\frac{1}{\lambda}\Pi
\mathbb{V}[\widetilde{w}^{(k-1)}(\lambda)-\widetilde{w}^{(k-1)}(0)],$$
where $$\widetilde{w}^{(1)}(0)=-R_0\mathbb{V}\varphi(u),$$

$$(\widetilde{w}^{(1)}(\lambda))'_{\lambda}|_{\lambda=0}=-R_0^2\mathbb{V}\Pi
\varphi(u),$$

$$\widetilde{w}^{(k)}(0)=R_0\mathbb{L}{u}^{(k-1)}(0)+R_0\mathbb{V}\widetilde{w}^{(k-1)}(0)+$$
$$\Pi
\mathbb{V}(\widetilde{w}^{(k-1)}(\lambda))'_{\lambda}|_{\lambda=0},$$

$$(\widetilde{w}^{(k)}(\lambda))'_{\lambda}|_{\lambda=0}=R^2_0\mathbb{L}{u}^{(k-1)}(0)+R_0^2Q_{1}\widetilde{w}^{(k-1)}(0)+$$
$$
R_0\mathbb{V}(\widetilde{w}^{(k-1)}(\lambda))'_{\lambda}|_{\lambda=0}.$$
\end{lem}

{\it Proof.}
$$\widetilde{w}^{(1)}(\lambda)=\int_0^{\infty}e^{-\lambda s}
{w}^{(1)}(s)ds=\int_0^{\infty}e^{-\lambda s} [e^{Qs}-\Pi]ds\cdot
w^{(1)}(0)=$$ $$=(\lambda-\Pi+(R_0+\Pi)^{-1})^{-1}[-\mathbb{V}
\varphi(u)],$$ where the correlation for the resolvent was found in
{\rm \cite{KorTur}}.

$$\widetilde{w}^{(1)}(0)=-R_0\mathbb{V}\varphi(u),$$
$$(\widetilde{w}^{(1)}(\lambda))'_{\lambda}|_{\lambda=0}=\lim\limits_{\lambda\to 0}\frac{R(\lambda)-R_0}{\lambda}[-\mathbb{V}
\varphi(u)]=-R_0^2\mathbb{V} \varphi(u).$$

For the next terms we have:
$$\widetilde{w}^{(k)}(\lambda)=(\lambda-\Pi+(R_0+\Pi)^{-1})^{-1}\mathbb{L}{u}^{(k-1)}(0)+$$
$$(\lambda-\Pi+(R_0+\Pi)^{-1})^{-1}\mathbb{V}
\widetilde{w}^{(k-1)}(\lambda)+\frac{1}{\lambda}\Pi
\mathbb{V}[\widetilde{w}^{(k-1)}(\lambda)-\widetilde{w}^{(k-1)}(0)],$$
here the last term was found using the following correlation:
$$\int_0^{\infty}e^{-\lambda s }\int_s^{\infty}
\mathbb{V}{w}^{(k-1)}(\theta)d\theta
ds=\int_0^{\infty}\int_0^{\theta}e^{-\lambda s
}\mathbb{V}{w}^{(k-1)}(\theta)ds d\theta=$$
$$\int_0^{\infty}\left(-\frac{1}{\lambda}\right)(e^{-\lambda\theta}
-1)\mathbb{V}{w}^{(k-1)}(\theta)d\theta =\frac{1}{\lambda}
\mathbb{V}[\widetilde{w}^{(k-1)}(\lambda)-\widetilde{w}^{(k-1)}(0)].$$

So,
$$\widetilde{w}^{(k)}(0)=R_0\mathbb{L}{u}^{(k-1)}(0)+R_0\mathbb{V}\widetilde{w}^{(k-1)}(0)+\Pi
\mathbb{V}(\widetilde{w}^{(k-1)}(\lambda))'_{\lambda}|_{\lambda=0},$$

$$(\widetilde{w}^{(k)}(\lambda))'_{\lambda}|_{\lambda=0}=R^2_0\mathbb{L}{u}^{(k-1)}(0)+R_0^2Q_{1}\widetilde{w}^{(k-1)}(0)+$$
$$R_0\mathbb{V}(\widetilde{w}^{(k-1)}(\lambda))'_{\lambda}|_{\lambda=0}-\lim\limits_{\lambda\to 0}\left\{\frac{1}{\lambda^2}\Pi
\mathbb{V}[\widetilde{w}^{(k-1)}(\lambda)-\widetilde{w}^{(k-1)}(0)]-\right.$$
$$\left.\frac{1}{\lambda}\Pi
\mathbb{V}(\widetilde{w}^{(k-1)}(\lambda))'_{\lambda}\right\},$$
where the last limit tends to 0.

Lemma is proved.

So, the obvious view of the initial condition for the $c^{(k)}(t)$
is:
$$c^{(k)}(0)=\mathbb{V}\widetilde{w}^{(k-1)}(0).$$

Theorem is proved.

\section{\bf Estimate of the remainder}
Let function $\varphi(u)$ in the definition of the functional
$\Phi_t^{\varepsilon}$ belongs to Banach space of twice
continuously differentiable by $u$ functions $C^2(\mathbb{R}^d)$.

Let us write (6) in the view
$$\tilde{\Phi}^{\varepsilon}(t)={\Phi}^{\varepsilon}(t)-{\Phi}^{\varepsilon}_2(t)\eqno(13)$$
where ${\Phi}^{\varepsilon}_2(t)=
{u}^{(0)}(t)+\varepsilon({u}^{(1)}(t)+{w}^{(1)}(t))+\varepsilon^2({u}^{(2)}(t)+
{w}^{(2)}(t)),$ and the explicit view of the functions
${u}^{(i)}(t), {w}^{(j)}(t), i=\overline{0,2}, j=1,2$ is given in
Theorem 2.1.

By theorem 3.2.1 from {\rm \cite{KorTur}} in Banach space
$C^2(R^d\times E)$ for the generator of Markovian evolution
$L^{\varepsilon}=\varepsilon^{-1}Q+\mathbb{V},$ exists bounded
inverse operator $(L^{\varepsilon})^{-1}.$

Let us substitute the function (13) into equation (6):
$$\frac{d}{dt}\tilde{\Phi}^{\varepsilon}-L^{\varepsilon}\tilde{\Phi}^{\varepsilon}=
\frac{d}{dt}{\Phi}^{\varepsilon}_2-L^{\varepsilon}{\Phi}^{\varepsilon}_2:=\varepsilon
\theta^{\varepsilon}. \eqno(14)$$

Here $\varepsilon
\theta^{\varepsilon}=\varepsilon[\frac{d}{dt}{u}^{(1)}-\varepsilon
\mathbb{V}({u}^{(2)}+{w}^{(2)})].$

The initial condition has the order $\varepsilon$, so we may write
it in the view:
$$\tilde{\Phi}^{\varepsilon}(0)=\varepsilon
\tilde{\Phi}^{\varepsilon}(0).$$

Let
$L_t^{\varepsilon}\varphi(u)=E[\varphi(u^{\varepsilon}(t))|u^{\varepsilon}(0)=u,
\ae^{\varepsilon}(0)=x]$ be the semigroup corresponding to the
operator $L^{\varepsilon}.$

\begin{thm} The following estimate is true for the remainder (13) of the
solution of equation (6):
$$||\tilde{\Phi}^{\varepsilon}(t)||\leq
 \varepsilon
||\tilde{\Phi}^{\varepsilon}(0)|| \exp\{\varepsilon
L||\theta^{\varepsilon}||\}, $$ where
$L\geq2||(L^{\varepsilon})^{-1}||.$
\end{thm}

{\it Proof:} The solution of equation (14) is:
$$\tilde{\Phi}^{\varepsilon}(t)=\varepsilon[L_t^{\varepsilon}
\tilde{\Phi}^{\varepsilon}(0)+ \int_0^t
L_{t-s}^{\varepsilon}\theta^{\varepsilon}(s)ds].$$

For the semigroup we have
$L_t^{\varepsilon}=I+L^{\varepsilon}\int_0^tL_s^{\varepsilon}ds,$
so
$\int_0^tL_s^{\varepsilon}ds=(L^{\varepsilon})^{-1}(L_t^{\varepsilon}-I).$

Using Gronwell-Bellman inequality {\rm \cite{BainSim}}, we receive
$$||\tilde{\Phi}^{\varepsilon}(t)||\leq \varepsilon L_t^{\varepsilon}
||\tilde{\Phi}^{\varepsilon}(0)||\exp\{\varepsilon\int_0^t
L_s^{\varepsilon}\theta^{\varepsilon}(t-s)ds\}\leq \varepsilon
L_t^{\varepsilon} ||\tilde{\Phi}^{\varepsilon}(0)||
\exp\{\varepsilon L||\theta^{\varepsilon}||\}, $$ where
$L\geq2||(L^{\varepsilon})^{-1}||.$

Theorem is proved.

\begin{rem} For the remainder of asymptotic expansion (5) of the view
$$\tilde{\Phi}^{\varepsilon}_{N+1}(t):={\Phi}^{\varepsilon}(t)-{\Phi}^{\varepsilon}_{N+1}(t),$$
where ${\Phi}^{\varepsilon}_{N+1}(t)=
{u}^{(0)}(t)+\sum_{k=1}^{N+1}\varepsilon^k({u}^{(k)}(t)+{w}^{(k)}(t))$
we have analogical estimate:
$$||\tilde{\Phi}^{\varepsilon}_{N+1}(t)||\leq
 \varepsilon^N
||\tilde{\Phi}^{\varepsilon}(0)|| \exp\{\varepsilon^N
L||\theta^{\varepsilon}_N||\},$$ where
$\frac{d}{dt}\Phi^{\varepsilon}_{N+1}-L^{\varepsilon}\Phi^{\varepsilon}_{N+1}:=\varepsilon^N
\theta^{\varepsilon}_N.$
\end{rem}

{\it Acknowledgements.}The author thanks Acad. V.S.Koroliuk for
the formulation of the problem studied. Acknowledgements also to
the Institute of Applied Mathematics, University of Bonn for the
hospitality and financial support by DFG  project 436 UKR
113/70/0-1.\par


\begin{thebibliography}{99}

\bibitem{BainSim} Bainov D., Simeonov P. {\it Integral inequalities
and applications,} Kluver Acad. Publ., Dordrecht, (1992), 316p.

\bibitem{Kor} Korolyuk V.S. {\it Boundary layer in asymptotic
analysis for random walks,} Theory of Stochastic Processes
\underline{1-2}, 25-36 (1998).

\bibitem{KorBor} Koroljuk V.S., Borovskikh Ju.V. {\it Analytic problems of
asymptotics of probabilistyc distributions,} Naukova umka, Kyiv,
(1981), 240 p. (in Russian).

\bibitem{KoKo} Korolyuk V.S., Korolyuk V.V. \textit{Stochastic Models
of Systems,} Kluwer Acad. Publ. (1999), 250p.

\bibitem{KoLim} Korolyuk V.S., Limnios N. \textit{Stochastic Systems in Merging Phase Space,}
World Scientific Publishers (2005), 330p.

\bibitem{KorPenTur} Koroljuk V.S., Penev I.P., Turbin A.F. {\it Asymptotic
expansion for the distribution of absorption time of Markov
chain,} Cybernetics \underline{4}, 133-135 (1973), (in Russian).

\bibitem{KorTur} Koroljuk V.S., Turbin A.F.
{\it Mathematical foundation of  state lumping of large systems,}
Kluver Acad. Press, Amsterdam, (1990), 280p.

\bibitem{Pin} Pinsky M. {\it Lectures on random evolutions,} World Scientific,
Singapore, (1991), 136 p.

\bibitem{Sam}  Samoilenko I.V. \textit{Asymptotic expansion for the
functional of markovian evolution in $R^d$ in the circuit of
diffusion approximation,} Journal of Applied Mathematics and
Stochastic Analysis 3, 247-258 (2005).

\bibitem{Skor} Skorokhod A.V., Hoppensteadt F.C., Salehi H.
\textit{Random Perturbation Methods with Applications in Science
and Engineering,} Springer (2002), 488p.

\bibitem{Tad} Tajiev A. {\it Asymptotic
expansion for the distribution of absorption time of semi-Markov
process,} Ukrainian Math. Journ. 9, 422-426 (1978). (in Russian)

\bibitem{VasBut} Vasiljeva A.B., Butuzov V.F. {\it Asymptotic
methods in the theory of singular perturbations,} Vyschaja shkola,
Moscow, (1990), 208 p. (in Russian).


\end{thebibliography}
\end{document}